\documentclass[twoside,12pt]{article}
\usepackage{amsmath,amssymb}
\date{}
\newtheorem{Lemma}{LEMMA}[section]
\newtheorem{Corollary}[Lemma]{COROLLARY}
\newtheorem{Theorem}[Lemma]{THEOREM}
\newtheorem{Proposition}[Lemma]{PROPOSITION}

\newcommand{\bnum}{\begin{enumerate}}
\newcommand{\enum}{\end{enumerate}}
\newcommand{\bi}{\begin{itemize}}
\newcommand{\ei}{\end{itemize}}
\newcommand{\btab}{\begin{tabular}}
\newcommand{\etab}{\end{tabular}}
\newcommand{\beq}{\begin{eqnarray*}}
\newcommand{\eeq}{\end{eqnarray*}}
\newcommand{\beqn}{\begin{eqnarray}}
\newcommand{\eeqn}{\end{eqnarray}}
\newcommand{\mat}[1]{\begin{matrix}#1\end{matrix}}

\newcommand{\bq}{\begin{equation}}
\newcommand{\eq}{\end{equation}}

\newcommand{\CL}{{\cal L}}

\def\phi{\varphi}
\def\epsilon{\varepsilon}

\newcommand{\kasten}{\vbox{\hrule height 8pt width 8.6pt depth -7.4pt
    \hbox{\vrule width 0.6pt height 7.4pt
    \kern 7.4pt \vrule width 0.6pt height 7.4pt}
    \hrule height 0.6pt width 8.6pt}}
\newcommand{\ok}{\hfill\kasten}
\newcommand{\bpf}{\begin{Proof}}
\newcommand{\epf}{\ok\end{Proof}\bigskip\noindent}
\newcommand{\bthm}{\begin{Theorem}}
\newcommand{\ethm}{\end{Theorem}}
\newcommand{\ble}{\begin{Lemma}}
\newcommand{\ele}{\end{Lemma}}
\newcommand{\bprop}{\begin{Proposition}}
\newcommand{\eprop}{\end{Proposition}}
\newcommand{\bcor}{\begin{Corollary}}
\newcommand{\ecor}{\end{Corollary}}
\begin{document}
\title{Compactness of the automorphism group of a topological parallelism on real projective 3-space: The disconnected case}

\author{Rainer L\"owen}

\maketitle
\thispagestyle{empty}

\begin{abstract}
We prove that the automorphism group of a topological parallelism on real projective 3-space is compact. 
This settles a conjecture stated in \cite{zush}, where it was proved that at least the connected component of the identity is compact.
 
\end{abstract}

\section{Introduction}

Clifford parallelism on real projective 3-space $\mathop {\rm PG}(3,\Bbb R)$ is the classical example of a 
topological parallelism. In general, a topological parallelism is an equivalence relation on the set of lines such that 
every equivalence class $\cal C$  is a spread (that is, its elements form a simple cover of the point set) and such that 
some continuity condition is satisfied. Examples of topological parallelisms abound, see, e.g., \cite{gl}, \cite{nonreg}, 
\cite{coll}, \cite{reg}, and further references given there. 
The group $\Phi = \mathop{\rm Aut} \Pi$ of automorphisms of a topological parallelism is a closed subgroup of 
$\mathop {\rm PGL}(4,\Bbb R)$. In the case of Clifford parallelism, $\Phi$ is the compact 6-dimensional group 
$\mathop{\rm PSO}(4,\Bbb R)$, see Corollary \ref{cliff} below. According to \cite{dim4}, Clifford parallelism is the only 
one with $\dim \Phi \ge 4$. Here we prove:

\bthm\label{main} 
Let $\Pi$ be a topological parallelism on $\mathop {\rm PG}(3,\Bbb R)$. Then the automorphism group $\Phi = \mathop{\rm Aut} \Pi$  
is compact and hence is conjugate to a closed subgroup of $\mathop{\rm PO}(4,\Bbb R)$.
\ethm

Correcting a statement in \cite{zush} we remark here that there are infinitely many conjugacy classes of 1-tori in 
$\mathop{\rm PO}(4,\Bbb R)$. The theorem will be proved in Sections 3 and 4 without using the result of \cite{zush}. Some
arguments are close to those used in \cite {zush} and will be treated somewhat sketchily or even be omitted. Other parts of the 
proof require new ideas. 

As an application of Theorem \ref{main}, we obtain a quick proof of the result of Betten and Riesinger 
\cite{cliff} on the automorphism 
group of Clifford parallelism:

\bcor\label{cliff}
{\rm \cite{cliff}} The automorphism group $\Delta$ of Clifford parallelism is the compact, connected 6-dimensional group 
$\mathop{\rm PSO}(4,\Bbb R)$.
\ecor

\bpf The group $\Sigma = \mathop{\rm PSO}(4,\Bbb R)$ is a direct product of two subgroups $\Lambda$ and $\Psi$ induced 
by the left and right multiplications by quaternions of norm one, respectively. Clifford parallelism is most easily 
defined as the set of all $\Lambda$-orbits in the line space of $\mathop {\rm PG}(3,\Bbb R)$. These orbits are permuted 
by the other factor $\Psi$ because $\Psi$ centralizes $\Lambda$. Thus $\Delta$  contains $\Sigma$. On the other hand, by
Theorem \ref{main}, $\Delta$ is contained in a maximal compact subgroup of $\mathop {\rm PGL}(4,\Bbb R)$, 
so up to conjugacy,
$\Delta$ is contained in $\mathop{\rm PO}(4,\Bbb R)$. Thus if $\Delta$ is disconnected, then it 
contains an involution $\gamma$
with determinant $-1$. 
Now $\gamma$ can be represented by a diagonal matrix with diagonal entries $-1, 1,1,1$, and it follows 
easily from Lemma \ref{fixed points} below that such an involution cannot leave any 
topological parallelism invariant. Thus we have equality $\Delta = \Sigma$.
\epf

\section{Preliminaries}

We consider real projective 3-space $\mathop {\rm PG}(3,\Bbb R)$ with its point space $P$ and line space $\CL$. 
Elements of $P$ and $\CL$ are the 1- and 2-dimensional subspaces of the vector space $\Bbb R^4$, respectively. 
The pencil of all lines passing through a point $p$ is denoted $\CL_p$.
The topologies of $P$ and $\CL$ are  the usual topologies of the Grassmann manifolds. 
This means that a sequence $X_n$ of subspaces converges to a subspace $X$ if and only if the spaces
$X_n$ have bases converging to a basis of $X$.
See, e.g.,  \cite{Kuehne} for the
continuity properties of geometric operations in this topological projective space. We shall mainly use the facts that $P$ and $\CL$ 
are compact and that the line $L = p \vee q$ joining 
two distinct points $p,q$ depends continuously on the pair $(p,q)$.
\\

A \it spread \rm on $\mathop {\rm PG}(3,\Bbb R)$ is a set of lines $\cal C \subseteq L$ such that every point 
is incident with exactly 
one of them. A spread defines an affine translation plane whose points are the elements of $\Bbb R^4$ and whose lines are 
the cosets of the vector spaces $L \in \cal C$. This translation plane is topological (i.e., joining points and 
intersecting lines are continuous operations) if and only if the spread $\cal C$ is a compact subset of $\CL$ 
(and then it is homeomorphic to the 2-sphere $\Bbb S_2$); see 
\cite{CPP}, 64.4.

A \it parallelism \rm  $\Pi$ on $\mathop {\rm PG}(3,\Bbb R)$ is a set of spreads such that every line is contained in exactly one 
of them. Thus $\Pi$ is an equivalence relation on $\CL$ whose classes are spreads. There are several ways of looking at 
a parallelism. Accordingly, the symbol $\Pi$ shall be used with different meanings. 
If considered as a relation on lines, $\Pi \subseteq \CL\times \CL$ is a set of pairs of lines. If considered as a set of spreads, 
$\Pi$ becomes a set of subsets of $\CL$. Finally, $\Pi$ may be considered as a map $\Pi: P \times \CL \to \CL$ sending a pair 
$(p,L)$ to the unique line $\Pi(p,L)$ parallel (i.e., equivalent) to $L$ that contains $p$, or as a map $\CL \to \Pi$ sending a line $L$ 
to the class $\Pi(L)$ containing $L$. We trust that this usage will not create any confusion.

A parallelism $\Pi$ is said to be \it topological \rm if the map $\Pi: P \times \CL \to \CL$ is continuous. It is not hard to see 
that this is equivalent to the condition that $\Pi \subseteq \CL \times \CL$ is a closed and, hence, compact subset, and that 
all classes of a topological parallelism are compact spreads. If we presuppose that the classes are compact spreads then $\Pi$
is topological precisely if $\Pi$ is a compact subset of the hyperspace $H(\CL)$ of all compact subsets of $\CL$, endowed with the 
Hausdorff topology or, equivalently, if the map $\Pi: \CL \to H(\CL)$ is continuous. Note that the restriction of this map to a line pencil $\CL_p$
is a homeomorphism from $\CL_p$ (which is a real projective plane) onto $\Pi \subseteq H(\CL)$. 

The automorphism group $\Phi = \mathop{\rm Aut}\Pi$ of a parallelism $\Pi$ is defined as the subgoup of 
$\mathop{\rm PGL}(4,\Bbb R)$ 
consisting of all collineations of $\mathop {\rm PG}(3,\Bbb R)$ that leave $\Pi$ invariant. 
Suppose that $\Pi$ is topological. 
From the compactness of $\Pi \subseteq H(\CL)$ one sees immediately that then $\mathop{\rm Aut}\Pi$ is 
a closed subgroup of the Lie group 
$\mathop{\rm PGL}(4,\Bbb R)$ and, hence, is again a Lie goup.

\section{Foundations of the proof}

The proof for the connected case given in \cite{zush} rests on the theorem of Malcev-Iwasawa, which implies that a 
non-compact Lie group contains some closed but non-compact one-parameter subgroup. This reduced the proof to the examination 
of all possible one-parameter groups. Here, we use instead the following theorem of Djocovi\'c: 

\bthm\label{djoc}
{\rm (\cite{djoc}, Theorem 3)} Let $G$ be a connected, locally compact group and let $H\le G$ be a closed subgroup. 
Suppose that every cyclic subgroup generated by an element of $H$ has compact closure (in $H$ or, equivalently, in $G$). 
Then $H$ is compact. 
\ethm  

We use this in order to reduce the proof of Theorem \ref{main} to the consideration of cyclic groups of automorphisms, as follows. 
Our group $G$ will be the connected component $\mathop{\rm PGL}(4,\Bbb R)^1$. It has index 2 in $\mathop{\rm PGL}(4,\Bbb R)$,
hence it suffices to show that $H = \Phi \cap \mathop{\rm PGL}(4,\Bbb R)^1$ is compact. This will follow from Theorem \ref{djoc} 
if we show that every cyclic subgroup of $\Phi$ has compact closure in $\mathop{\rm PGL}(4,\Bbb R)$. 

It suffices to consider one cyclic group from every conjugacy class. In order to obtain these representatives, we use the 
Jordan normal form of real $4 \times 4$ matrices. We also did this in \cite{zush}, but there we applied this method to 
generators of one-parameter groups rather than generators of cyclic groups. Consequently, we may now replace a matrix by 
any nonzero scalar multiple as we are dealing with the projective group, whereas in \cite {zush} we were allowed to 
\it add \rm an arbitrary scalar matrix. This leads to a different list of candidates. In compiling this list, we observe 
that all eigenvalues will be nonzero. We may replace a matrix by its square without loss, 
hence we may assume that real eigenvalues are positive.
In the list which follows, the cases labelled `a' are the complex $2 \times 2$ matrices, label `b'  indicates 
the presence of both real and complex eigenvalues, and label `c' is attached to matrices all of whose eigenvalues are real.
The case distinction is not identical to the one used in \cite{zush}. We allow complex eigenvalues to be real 
eigenvalues with a 2-dimensional eigenspace, and in this case we use a complex coordinate on that eigenspace. 
We shall write $\gamma$ for a real $4 \times 4$ matrix and sometimes also for the projective collineation induced by that matrix. 
The cyclic group generated by $\gamma$ is denoted $\langle \gamma \rangle$. \\

\bf Case (a1):  $\gamma$  \it is a complex diagonal $2 \times 2$ matrix with (possibly real) eigenvalues 
$a, c \in \Bbb C\setminus\{0\}$: \rm

      $$\gamma = \left(
      \begin{matrix} a & \cr
                      & c \cr 
                      \end{matrix}
         \right)  .            $$
Up to multiplication by a real number, we have that $\vert a \vert = 1$. If also $\vert c \vert = 1$, then $\langle \gamma \rangle$
is contained in the 2-torus of all complex diagonal matrices with diagonal entries of norm one and, hence, has compact closure. 

\bf Case (a2): $\gamma$ \it is a complex Jordan block with (possibly real) eigenvalue $a \in \Bbb C\setminus \{0\}$. \rm             

Multiplying by $b = \vert a \vert^{-1} \in \Bbb R$ and renaming, we get $\vert a \vert = 1$ and 
 
 $$\gamma = \left(
      \begin{matrix} a & b \cr
                      & a \cr 
                      \end{matrix}
         \right), \quad  
               \gamma^n = \left(
      \begin{matrix} a^n & na^{n-1}b \cr
                      & a^n \cr 
                      \end{matrix}
         \right)  .   $$
(To compute the powers of $\gamma$, write $\gamma$ as a sum of a scalar matrix plus a nilpotent matrix and use the 
binomial theorem. This works in more complicated cases as well, cf. case (c1) below.)         
   
\bf Case (b1): \it $\gamma$ has one eigenvalue $a \in \Bbb C \setminus \{0\}$ and two real ones 
$r,s \in \Bbb R\setminus \{0\}$, the latter with Jordan blocks of size 1. \rm

If $a$ happens to be real, we assume that it has an eigenspace of real dimension 2. 
As before, we assume $\vert a \vert = 1$ and obtain

 $$\gamma = \left(
      \begin{matrix} a & &\cr
                      & r & \cr
                     && s \cr 
                      \end{matrix}
         \right)  .            $$
         
\bf Case (b2): \it $\gamma$ has one eigenvalue $a \in \Bbb C \setminus \{0\}$ and one real eigenvalue $r$, the 
latter with a Jordan block of size 2.    \rm

Again $a$ is assumed to have a 2-dimensional eigenspace if it is real. Dividing by $\vert a \vert$ and renaming, we get

     $$\gamma = \left(
      \begin{matrix} a & &\cr
                      & r & s \cr
                     && r \cr 
                      \end{matrix}
         \right)  , 
          \quad {\rm where}\  \vert a \vert = 1 \quad {\rm and} \quad r,s \in \Bbb R \setminus \{ 0 \}.        $$

\bf Case (c1): \it $\gamma$ has one real eigenvalue $s \ne 0$ with a single Jordan block of size 4. \rm         

Multiplying by $r = s^{-1}\in \Bbb R \setminus \{ 0 \}$, we get 

     $$\gamma = \left(
      \begin{matrix} 1 & r &  & \cr
                      & 1 & r &  \cr
                     && 1 & r \cr 
                     &&& 1 \cr
                      \end{matrix}
         \right), \quad 
         \gamma^n = \left(
      \begin{matrix} 1 & nr & \binom{n}{2} r^2 & \binom{n}{3} r^3\cr
                      & 1 & nr & \binom{n}{2} r^2 \cr
                     && 1 & nr \cr 
                     &&& 1 \cr
                      \end{matrix}
         \right)  .            $$
         
\bf Case (c2): \it $\gamma$ has two real eigenvalues with Jordan blocks of size 3 and 1, respectively. \rm

After dividing by one of the eigenvalues, we have

 $$\gamma^n = \left(
      \begin{matrix} 1 & nr & \binom{n}{2} r^2 & \cr
                      & 1 & nr &  \cr
                     && 1 &  \cr 
                     &&& s^n \cr
                      \end{matrix}
         \right) , \quad r,s \in \Bbb R \setminus \{0\}.            $$

\bf Case (c3): \it $\gamma$ has two distinct real eigenvalues, both with one Jordan block of size 2. \rm    

(The corresponding case with two equal eigenvalues is covered by case (a2), up to a change of coordinates.)
By passing to $\gamma^2$ we obtain that the eigenvalues are positive. After dividing by one of the eigenvalues, we have

      $$\gamma = \left(
      \begin{matrix} r & s & & \cr
                      & r &  &  \cr
                     && 1 & s \cr 
                     &&  & 1 \cr
                      \end{matrix}
         \right) , \quad
\gamma^n = \left(
      \begin{matrix} r^n & nr^{n-1}s & & \cr
                      & r^n &  &  \cr
                     && 1 & ns \cr 
                     &&  & 1 \cr
                      \end{matrix}
         \right) , \quad 0 < r,s \in \Bbb R, r \ne 1.            $$
         
By inverting $\gamma$, if necessary, we can enforce $0 < r <1$.         
       
\bf Case (c4): \it $\gamma$ has three real eigenvalues  with Jordan blocks of size 2, 1, 1, respectively. \rm

After dividing by one of the eigenvalues, we have

 $$\gamma = \left(
      \begin{matrix} 1 & t & & \cr
                      &1 &  &  \cr
                     && r &  \cr 
                     &&  & s \cr
                      \end{matrix}
         \right), \quad 0< r, s, t \in \Bbb R.            $$
We may assume that $r\ne s$, the case $r=s$ being covered by case (b2).\\         

\bf Case (c5): \it $\gamma$ is diagonal with four distinct eigenvalues $1,r,s,t$. \rm
(Again if two eigenvalues are equal, we have case (b1).)\\

We close this section by providing some tools that will be useful in the main part of the proof and possibly elsewhere. \\

\ble \label{disjoint}
Given a topological parallelism $\Pi$ and $m+1$ lines $M, L_1,...,L_m$, there is always a line $N$ parallel to $M$ and not 
meeting any of the $L_i$.
\ele

\bpf Lines $L_j$ that happen to be parallel to $M$ are easily dealt with. Suppose that this is not the case for $L_i$.
Then the set of parallels to $M$ that meet $L_i$ is the image of the point set of $L_i$ under the injective map $p \mapsto \Pi(p,M)$, 
hence like $L_i$ it is a topological circle. The union of a finite number of circles must be a proper subset of the 2-sphere $\Pi(L)$.
All lines $N$ in the complement satisfy our assertion.\epf 

\ble\label{existence}
Let $\Pi$ be a topological parallelism and let $K,L$ be two lines with trivial intersection. 
Then there exist two distinct parallel lines
$M,N$ meeting both $K$ and $L$.
\ele

\bpf 
The point sets $[K]$ and $[L]$ are topological circles. If the map $(p,q) \to \Pi(p\vee q)$ from the torus $[K] \times [L]$ into 
$\Pi\subseteq H(\CL)$ is injective, then it is also surjective by domain invariance, and we have a 
homeomorphism between a torus and a projective 
plane. This is impossible, hence there are two pairs $(p,q)$ and $(p', q')$ whose joining lines are parallel.
\epf

The next lemma is a variation of \cite{zush}, 3.8.

\ble\label{fixed points} Suppose that $\gamma$ leaves $\Pi$ invariant and fixes 
two points $p,q \in P$. Then the actions induced by $\gamma$
on the factor spaces $\Bbb R^4/p$ and $\Bbb R^4 / q$ are projectively equivalent.
\ele

\bpf The points of the projective space associated with the vector space $\Bbb R^4/p$ 
are in $1-1$ correspondence to the lines of $\mathop{\rm PG}(3,\Bbb R)$ passing through $p$.
The map $L \to \Pi(q,L)$ is  a $\gamma$-equivariant bijection of line pencils $\CL_p \to \CL_q$. 
\epf

\section{Completion of the proof}

We look at the matrices $\gamma$ listed in the previous section one by one. In each case we shall show that either the cyclic group
$\langle \gamma \rangle$ has compact closure in $\mathop{\rm GL}(4,\Bbb R)$  or $\gamma$ cannot belong to the 
automorphism group $\Phi= \mathop{\rm Aut} \Pi$ of a topological parallelism. Types (ai) and (bi) are the ones that require new ideas.

\bprop\label{a1}
If $\Phi$ contains an element $\gamma$ of type (a1), then both eigenvalues have norm 1, and $\langle \gamma \rangle$ has compact closure.
\eprop

\bpf 1) We write the second eigenvalue as $c = rb$ with $0 < r\in \Bbb R$ and $\vert b \vert = 1$. If $r = 1$, then $\gamma$ 
belongs to the 2-torus $T$ of unitary diagonal matrices. If $r \ne 1$, we may assume that $r > 1$. Consider the diagonal matrix 
$\tau \in T$ with eigenvalues  $a, b$. If this has infinite order, then its powers must accumulate at some element of $T$. 
Hence in any case, there is a sequence $n(k)$ of positive integers tending to infinity and such that the powers 
$\tau^{n(k)}$ converge to the unit matrix. We set $\delta_k := \gamma^{n(k)}$. 

2) We write elements of $\Bbb C^2$ as column vectors and consider a line 
    $$M = \left \langle \left(\begin{matrix} x \cr 0 \end{matrix} \right ), 
               \left( \begin{matrix} 0\cr y \end{matrix} \right) \right\rangle$$ 
meeting both factors of $\Bbb C \times \Bbb C$ nontrivially.  The line $M^{\delta_k}$ has basis vectors 
$$\left ( \begin{matrix} a^{n(k)}x \cr 0 \end{matrix} \right ) , \left ( \begin{matrix} 0 \cr b^{n(k)}y \end{matrix} \right )$$ 
and hence converges to $M$ as $k \to \infty$.

3) On the other hand, a line $N$ disjoint from $\Bbb C \times 0$ may be described using a real $2 \times 2$ matrix $A$; complex 
numbers are thought of as column vectors in $\Bbb R^2$ when they are multiplied by the matrix $A$: 
    $$N = \left \{ \left( \begin {matrix} Ay \cr y \end{matrix} \right )  ;  \ y \in \Bbb C\right\}.$$ 
The images $N^{\delta_k}$ converge to $0 \times \Bbb C$ because $y$ is arbitrary and 
   $$ r^{-n(k)} \left ( \begin{matrix} Ay \cr y \end{matrix} \right )^{\delta_k} 
        \to \left (\begin {matrix} 0 \cr y \end{matrix} \right ) .$$   
4) Now the line $M$ considered above cannot be parallel to the line $0 \times \Bbb C$ since the two lines have a point in common. 
However by Lemma \ref{disjoint}, there is a line $N$ parallel to $M$ and disjoint 
from $\Bbb C \times 0$. 
Then the $\delta_k$-images of $M$ and $N$ 
are also parallel, and passing to the limit we obtain that $M$ is parallel to $0 \times \Bbb C$, a contradiction.
\epf

\bprop \label{a2}
The group $\Phi$ does not contain any nontrivial element $\gamma$ of type (a2). 
\eprop

\bpf 1) As in the previous proof, we take a sequence $n(k)$ of positive integers such that $n(k) \to \infty$ and $a^{n(k)} \to 1$,
and we set $\delta_k := \gamma^{n(k)}$. 

2) For $x,y \in \Bbb C, y \ne 0$ we have
   $$(n(k)b)^{-1}\left ( \mat{x \cr ay} \right )^{\delta_k} =
      \left (\mat{(n(k)b)^{-1}a^{n(k)}x + a^{n(k)}y \cr (n(k)b)^{-1}a^{n(k)+1}y} \right ) \to \left ( \mat{y \cr 0}\right ).$$

3) Now consider the line 
    $$L = \left \langle \left ( \mat{1 \cr 0} \right ), \left (\mat{0 \cr a}\right ) \right \rangle.$$
By step (2), if $L^{\delta_k}$ accumulates at any line $X$ then $X$ contains    $\left ( \mat{1 \cr 0}\right )$. 
Moreover, $X$ contains the limit of the sequence 
    $$\left ( \mat{-n(k)b \cr a}\right ) ^{\delta_k} = \left( \mat{0 \cr a^{n(k)+1}} \right ),$$
which is equal to $\left ( \mat{0\cr a}\right )$. Thus we see that $L^{\delta_k} \to L$ as $k \to \infty$. 
 
4) By Lemma \ref{disjoint}, there is a line $M$ parallel to $L$ that has trivial intersection with $W : = \Bbb C \times 0$.
Then $M$ contains two vectors whose second components are $\Bbb R$-linearly independent, hence $M^{\delta_k} \to W$ by
step (2).
Applying $\delta_k$ to the parallel pair $(L,M)$ and passing to the limit, 
we obtain that $L$ is parallel to $W$, a contradiction because $L$ and $W$ have
a point in common.   $\phantom x$
\epf

\bprop\label{b1} If \ $\Phi$ contains an element $\gamma$ of type (b1), then the eigenvalues $a, r, s$ have norm 1, and 
$\langle \gamma \rangle$ has compact closure.
\eprop

\bpf We consider the case that not both of $r,s$ have modulus 1 and aim for a contradiction.

1) Assume first that $r = s$. Taking the square of $\gamma$, we obtain $r > 0$, and passing to $\gamma^{-1}$, if necessary, we get $r>1$.
We use a sequence $n(k)$ of integers as in the previous proofs and also retain the definition of $\delta_k$. 
We denote the lines $\Bbb C \times 0$ and $0 \times \Bbb C$ by $W$ and $S$, respectively.

2) Since we assume that $r = s$, a line $L$ that meets both $S$ and $W$ has images $L^{\delta_k}$ converging to $L$ itself. 
On the other hand,  a line $M$ that intersects $S,W$ trivially may be described, using a matrix $A \in \mathop{\rm GL}(2,\Bbb R)$, as
    $$M = \left \{ \left ( \mat{x \cr Ax}\right ) ; \quad x \in \Bbb C\right \},$$
and for its images we have
    $$M^{\delta_k} = \left \{ \left (\mat{a^{n(k)}x \cr r^{n(k)}Ax}\right ) ; \quad x \in \Bbb C\right \} \to
     \left \{ \left (\mat{0 \cr Ax}\right ) ; \quad x \in \Bbb C\right \}  = S.$$
By Lemma \ref{disjoint},  the line $M$ can be chosen to be parallel to $L$, and passing to the 
limit we obtain that $L$ is parallel to $S$, a contradiction. 

3) Now suppose that $r \ne s$. We denote by $e_i$ the standard basis vector of $\Bbb R^4$ having coordinates 0 in position $j \ne i$ and coordinate 
1 in position $i$. The corresponding points will be denoted $p_i = \langle e_i\rangle$. Applying Lemma \ref {fixed points} to $(p,q) = (p_3,p_4)$, 
we see that $1 \notin \{r,s\}$ and also that 1 does not lie between $r$ and $s$. Thus we may assume that $1<r<s$. This leads to a 
contradiction as in the proof of \cite{zush}, 3.13. For the convenience of the reader, we sketch the argument: 
We let $K = p_2 \vee p_4$, $L = p_3 \vee p_4 = L^{\delta_k}$. We have $p_i^{\delta_k} \to p_i$ for $i \in \{2,4\}$, hence $K^{\delta_k} \to K$. 
We consider two hyperplanes $G = e_1^\perp$ and $H = e_4^\perp$. By Lemma \ref {disjoint} there is a line $M$ parallel to $K$ but not equal to $K$
and not intersecting the lines $G\cap H$ or $p_1 \vee p_2$. Thus $M$ may be written as $M = p \vee q$, where $p$ is a point incident 
with $G$ but not with $H$, and 
$q$ is a point incident with $H$ but not with $(p_1 \vee p_2)$.  Then $p^{\delta_k} \to p_4$ and $q^{\delta_k}\to p_3$, hence $M^{\delta_k} \to L$.
However, 
   $$M^{\delta_k} = \Pi(p,K)^{\delta_k} = \Pi(p^{\delta_k},K^{\delta_k}) \to \Pi(p_4,K) = K.$$
This contradicts   $M^{\delta_k} \to L$. 
\epf

\ble\label{b2} The group $\Phi$ does not contain any nontrivial element $\Gamma$ of type (b2).
\ele

\bpf  As usual, we may assume that $r>1$, and we have sequences $n(k)$ and $\delta_k$ with the usual properties. 
We consider the lines $W= \Bbb C \times 0$ and $S = 0\times \Bbb C$, which have no common point. If a line $L$ meets
both $W$ and $S$, then there are
nonzero complex numbers $x,y \in \Bbb C$ such that
   $$ L = L(x,y) = \left \langle \left  ( \mat{x \cr 0}\right ), \left ( \mat{0\cr y}\right ) \right \rangle.$$
It is easily seen that
   $$\left (\mat{x \cr 0}\right ) ^{\delta_k} \to \left (\mat{x \cr 0}\right ) \quad {\rm and} \quad
     \left\langle \left (\mat{0 \cr y}\right ) \right \rangle ^{\delta_k} \to \left \langle \left (\mat{0 \cr 1}
     \right )\right \rangle,$$
whence
      $$L(x,y)^{\delta_k} \to M(x) = \left \langle \left (\mat{x\cr 0}\right ), \left (\mat{0 \cr 1}\right )\right \rangle.$$   
According to Lemma \ref{existence}, there is a pair of distinct parallel lines $L(x,y)$, $L(x',y')$. 
Then $x\ne x'$, but $M(x)$ is parallel to $M(x')$ by continuity, a contradiction because these lines have a common point.     
\epf

In contrast to the cases with complex eigenvalues treated so far, the remaining cases with only real eigenvalues do 
not require completely new ideas but can be treated by adapting the proofs given in the connected case in \cite{zush}.
We shall describe the necessary changes and in one case we shall give an alternative proof. Some simplification 
over the proof in \cite {zush} is 
achieved by shifting subcases from the `c' set to the `a' or `b' sets.

\bprop\label{c}
If the group $\Phi$ contains an element $\gamma$ of type $(ci)$, $1 \le i \le5$, then $\gamma$ has order at most two.
\eprop

\bpf
The square of $\gamma$ has positive eigenvalues, hence it suffices to show that $\Phi$ does not contain a nontrivial 
element of type (ci) with positive eigenvalues. 

\it Type (c1). \rm This is the most subtle part of the whole proof. Therefore, we briefly give the argument, although it is almost identical 
to the one used in \cite{zush}, 3.5. The matrix $\gamma$ can be written as $1+R$, where $1$ denotes the unit matrix and 
    $$R = \left(\mat{0&r& &\cr &0&r&\cr &&0&r \cr &&&0}\right).$$
Using the binomial theorem for commuting matrices together with the relation $R^4 = 0$, we obtain 
$\gamma^n = 1+nR + {1\over 2}(n^2 - n)R^2 + {1\over 6}(n^3-3n^2+2n)R^3$ for positive integers $n$. Moreover,
we find $\gamma^{-n} = 1 - nR + {1\over 2} (n^2+n)R^2 - {1\over 6}(n^3-9n^2+2n)R^3$. The coefficients are polynomials in $n$ whose degree increases as 
we move away from the main diagonal of the matrix.  From this observation it follows that the dynamics of $\langle \gamma \rangle$ acting on 
$P$ and on $\CL$ is the same as for a nontrivial cyclic subgroup of a one-parameter group of type (c1) as considered in \cite{zush}.
This means that $p^{\gamma^{\pm n}}$ converges to $p_1 = \langle e_1\rangle$ as $n \to \infty$, for points $p \notin H = e_4^\perp$.
Moreover, if $L_n$ is a sequence of lines that contain $p_1$ and the sequence converges to a line $L$ not contained in $H$, then $L_n^{\gamma^{\pm n}}$ 
converges to $L_{12}$, where we set $L_{ij} = \langle e_i, e_j\rangle$. We say that $\gamma^{\pm n}$ continuously converges to the constant map 
with value $L_{12}$.

Now we let $s_n = p_3^{\gamma^{-n}}$ and $M_n = s_n \vee p_4$; then $s_n \to p_1$ and, hence, $M_n \to L_{14}$. The images 
$M_n^{\gamma^n} = p_3 \vee p_4^{\gamma^n}$ converge to $L_{13}$. By continuity of the parallelism, 
the parallel lines $\Pi(M_n^{\gamma^n},p_1)$ converge to 
$\Pi(L_{13},p_1) = L_{13}$. But on the other hand, $\Pi(M_n, p_1)$ contains $p_1$ and converges to $\Pi(L_{14},p_1) = L_{14}$, a 
line not contained in $H$. Thus by continuous convergence, $\Pi(M_n^{\gamma^n}, p_1) = \Pi (M_n,p_1)^{\gamma^n}$ 
converges to $L_{12}$, a contradiction.
Thus, type (c1) cannot occur.

\it Type (c2). \rm The argument from \cite{zush}, 3.7 goes through, but we give a variation. The line 
$L = \langle e_1,e_4\rangle$ is fixed by $\gamma$. There is a unique line $K$ parallel  to $L$ which is contained in the 
hyperplane $H = \langle e_1,e_2,e_3\rangle$ corresponding to the Jordan block of size 3, see \cite{CPP}, 64.10a.
Hence $K$ is fixed by $\gamma$ and must contain an eigenvector for the unique eigenvalue of $\gamma$ on $H$. 
Thus $e_1 \in K$, a contradiction.

\it Type (c3). \rm  The proof from \cite{zush} can be adapted, but we prefer to sketch a more straightforward 
(though perhaps less beautiful) argument.

Any line  which has trivial intersection with $W = \Bbb C \times 0$ can be described as
    $$L_A = \left \{ \left (\mat{Ay \cr y}\right ) ; \ y \in \Bbb C \right \}, $$
where $A$ is a real $2 \times 2$ matrix. A block diagonal matrix 
   $$\left (\mat {B & \cr & C} \right )$$
transforms $L_A$ into $L_{BAC^{-1}}$, and  an easy computation shows in view of $0<r<1$ that $L_A^{\gamma^n}$ 
converges to $L_0 = S = 0 \times \Bbb C$. Now by \ref{disjoint}, the line $K = \langle e_1,e_3\rangle$, which is 
fixed by $\gamma$,  is parallel to some line $L_A$. This implies that $K$ is parallel to $S$, a contradiction.   

\it Type (c4). \rm As in \cite{zush}, this is easily ruled out by applying Lemma \ref{fixed points} to the points 
$\langle e_1\rangle$ and $\langle e_3\rangle$.

\it Type (c5). \rm Both proofs of \cite{zush}, 3.13 can be used without change. Note that we need only consider the case where 
the eigenvalues are pairwise different. 
\epf

The proof of Theorem \ref{main} is now complete. \\

\noindent \bf Acknowledgement. \rm I am indebted to Harald L\"owe for valuable conversations and for drawing my 
attention to Djocovi\'c's result \ref{djoc}.


\bibliographystyle{plain}

\bigskip
\bigskip
\noindent Rainer L\"owen, Institut f\"ur Analysis und Algebra,
Technische Universit\"at Braunschweig,
Universit\"atsplatz 2,
D 38106 Braunschweig,
Germany

\end{document}